\documentstyle[11pt]{article}
\newlength{\spacing}
\newcommand{\doublespace}{\setlength{\baselineskip}{1.5\spacing}}
\newtheorem{thm}{Theorem}[section]

\newtheorem{lem}[thm]{Lemma}
\newtheorem{prop}[thm]{Proposition}
\newtheorem{defn}{Definition}

\def\lab{\label }

\def\rar{\to}

\def\inft{\infty}

\def\del{\delta}

\def\al{\alpha}

\def\ep{\epsilon}

\def\today{\ifcase\month\or
  January\or February\or March\or April\or May\or June\or
  July\or August\or September\or October\or November\or December\fi
  \space\number\day, \number\year}

\begin{document}
\begin{titlepage}
\begin{center}
{\bf Vertex Degree of Random Intersection Graph} \\
\vspace{0.20in} by \\
\vspace{0.2in} {Bhupendra Gupta \footnote{Corresponding Author.
email:gupta.bhupendra@gmail.com, bhupen@iiitdm.in}}\\
Faculty of Engineering and Sciences,\\
Indian Institute of Information Technology (DM)-Jabalpur, India.\\
\vspace{0.1in}
\end{center}
\vspace{0.2in}
\sloppy
\begin{center} {\bf Abstract} \end{center}

\begin{center} \parbox{4.8in}
{
A random intersection graph is constructed by independently
assigning a subset of a given set of objects $W,$ to each vertex of
the vertex set $V$ of a simple graph $G.$ There is an edge between
two vertices of $V,$ iff their respective subsets(in $W$,) have at
least one common element. The strong threshold for the connectivity
between any two arbitrary vertices of vertex set $V,$ is derived.
Also we determine the almost sure probability bounds
for the vertex degree of a typical vertex of graph $G.$
} \\
\vspace{0.4in}
\end{center}

\vspace{0.5in}
\vspace{0.5in} {\sl AMS 2000 subject classifications}:
\hspace*{0.5in} 05C80, 91D30.
\end{titlepage}
\doublespace
\section{Introduction\lab{s1}}
The development in this paper is in continuation of Dudley Stark
\cite{Dudley}, in which distribution of the degree of a typical
vertex is given. Dudley studied the model given in \cite{Singer}
defined as:\\ 
Let us consider a set $V$ with $n$ vertices and another set of
objects $W$ with $m$ objects. Define a bipartite graph $G^*(n,m,p)$
with independent vertex sets $V$ and $W.$ Edges between $v \in V$
and $w \in W$ exists independently with probability $p.$ The random
intersection graph $G(n,m,p)$ derived from $G^*(n,m,p)$ is defined
on the vertex set $V$ with vertices $v_1,v_2 \in V$ are adjacent if
and only if there exists some $w \in W$ such that both $v_1$ and
$v_2$ are adjacent to $w$ in $G^*(n,m,p).$ Also define $W_v$ be a
random subset of $W$ such that each element of $W_v$ is adjacent to
$v \in V.$ Any two vertices $v_1,v_2 \in V$ are adjacent if and only
if $W_{v_1} \cap W_{v_2} \neq \phi,$ and edge set $E(G)$ is define
as
\[E(G) = \{\{v_i,v_j\}: v_i,v_j \in V, \:W_{v_i} \cap W_{v_j} \neq \phi\}.\]
\section{Definitions and Supporting Results\lab{s2}}
\begin{lem}
Let $X$ be a random variable having binomial distribution with
parameters $n$ and $p>0,$ i.e., $X \sim Bi(n,p).$ Then
\begin{equation}
P[X \geq k] \leq \left(\frac{np}{k}\right)^k \exp(k-np), \qquad
k\geq np,\label{e1}
\end{equation}
and
\begin{equation}
P[X \leq k] \leq \left(\frac{np}{k}\right)^k \exp(k-np), \qquad 0 <
k\leq np.\label{e2}
\end{equation}
 \label{lemma1}
\end{lem}
For convenience, the bounds (\ref{e1}) and (\ref{e2}) can be
expressed as
\[\exp\left(npH\left(\frac{np}{k}\right)\right),\]
where
\[H(t) = \frac{1}{t}\log\:t + \frac{1}{t} -1,\quad 0<t<\inft,\]
and $H(\inft) = -1.$ Note that $H(t)<0,\: \forall \:\: t \in
(0,\infty]\setminus\{1\};\: H$ is increasing on $(0,1)$ and
decreasing on $(1,\infty).$ (See Lemma 1.2 on page 25, Penrose \cite{Penrose}.) \hfill$\Box$\\

Let $X = X(n,m,p)$ be the number of vertices of $V-\{v\}$ adjacent
in $G(n,m,p)$ to a vertex $v \in V,$ i.e., $X$ be the vertex degree
of a vertex $v \in V$ in $G(n,m,p).$ Then $X$ follows Bi$(n-1,
q_n),$ where $q_n$ is the probability of $v_1,v_2 \in V$ are adjacent.\\

\begin{prop}
Let $v_i,\:v_j\in V,\: i\neq j$ and $i,j = 1,2,\ldots,n.$ Then $v_i$
and $v_j$ are adjacent with probability $q_n,$ such that
\[q_n \sim mp^2,\]
for sufficiently small $p.$
\label{prop1}
\end{prop}
\textbf{Proof.} Consider
\begin{eqnarray*}
q_n & = & P[v_i, v_j \mbox{ are adjacent.}]\nonumber\\
& = & P[W_{v_i}\cap W_{v_j} \neq \phi]\nonumber\\
& = & P[\mid W_{v_i}\cap W_{v_j} \mid \geq 1]\nonumber\\
& = & 1- P[ w \notin (W_{v_i}\cap W_{v_j})]^m\nonumber\\
& = & 1-[1- P[ w \in (W_{v_i}\cap W_{v_j})]]^m\nonumber\\
& = & 1-[1- P[ w \in W_{v_i}, w \in W_{v_j})]]^m\nonumber\\
& = & 1-[1-p^2]^m
%
\end{eqnarray*}
Using the Taylor's series expansion up to second term, we get
\begin{equation}
q_n = mp^2 - \zeta_m, \label{q}
\end{equation}
where, $\zeta_m = \frac{m(m-1)}{2!}(1-c)^2p^4,$ and $c \in
(0,p^2).$

Now if we take $p$ is sufficiently small. Then
\begin{equation}
q_n \sim mp^2.
\end{equation}
This completes the proof. \hfill$\Box$\\

%
%
\begin{defn}
Let graphs $A$ and $B$ share the same vertices and the edge set of
$A$ is a subset of the edge set of $B$, we write $A \leq B.$ Also
let $\Theta$ be a property of a random graphs such that if $A \leq
B$ and $A \in \Theta,$ then $B \in \Theta.$ (Here $A \in \Theta$ is
used to denote that graph $A$ has property $\Theta$.) Then $\Theta$
is called an {\it upwards-closed property}. If $B \in \Theta$
implies $A \in \Theta$, then $\Theta$ is said to be a {\it
downwards-closed property}. If property $\Theta$ is upwards-closed
or downwards-closed, then $\Theta$ is called {\it monotone property.}\hfill$\Box$\\
\end{defn}

Fix a monotone property $\Theta.$ For any two functions $\delta,
\gamma: Z^+ \rar {R}^+$, we write $\delta \ll \gamma$ (resp. $\del
\gg \gamma$) if $\del(n)/\gamma(n) \rar 0,$ (resp.
$\gamma(n)/\delta(n) \rar 0$) as $n \rar \infty.$ We will write
$\delta$ for $\delta(n).$ Let $G_n(r)$, be any random graph on $n$
vertices.

\begin{defn}
  A function $\delta_{\Theta}: Z^+ \rar {R}^+$ is a {\it weak
  threshold}
  function for $\Theta$ if the following is true for every function
  $\delta: Z^+ \rar {R}^+,$
  \begin{itemize}
  \item if $\delta(n) \ll \delta_{\Theta}(n)$, then $P[G_n(\delta)
    \in \Theta] = 1-o(1),$ and
  \item if $\del(n) \gg \del_{\Theta}(n)$ then $P[G_n(\del) \in
    \Theta] = o(1)\ .$\hfill$\Box$\\
  \end{itemize}

  A function $\del_{\Theta}: Z^+ \rar {R}^+$ is a {\it strong
  threshold}
  function for $\Theta$ if the following is true for every fixed
  $\ep>0,$
 \begin{itemize}
  \item if $P[G_n(\del_{\Theta}-\ep) \in \Theta] = 1 - o(1),$ and
  \item if $P[G_n(\del_{\Theta}+\ep) \in \Theta] = o(1)\ .$\hfill$\Box$\\
  \end{itemize}
\end{defn}
\section{Main Results}\lab{s3}
We now derive the threshold probability for the connectivity of
graph $G(n,m,p).$

Let $\Theta$ be the connectivity of graph $G(n,m,p)$ and define  the
probability that a vertex $v \in V$ is connected to $w \in W,$ as
follows
\[p
:= p(\al) = \frac{1}{(mn^{\al})^{1/2}}\:\:.\]
\begin{thm}
Let $G(n,m,p)$ be the random intersection graph with $p(\al)=
\frac{1}{(mn^{\al})^{1/2}}$ and $v_i, v_j \in V(G)$ for $i \neq j =
1,2,\ldots,n.$ Then $p(2) = \frac{1}{m^{1/2}n},$ is a strong
threshold for the random intersection graph $G(n,m,p),$ i.e.,
\begin{eqnarray}
P[G(n,m,p(2-\ep)) \in \Theta] & = & 1-o(1),\nonumber\\
P[G(n,m,p(2+\ep)) \in \Theta] & = & o(1).\label{thres}
\end{eqnarray}
%
%
\label{conn}
\end{thm}
\textbf{Proof.} Let $v_i, v_j \in V(G)$ for $i \neq j,\: i,j =
1,2,\ldots,n.$ Then from (\ref{q}), we have
\begin{eqnarray}
P[\mbox{$v_i,v_j$ are adjacent}] & = & q_n\nonumber\\
& \leq & mp^2.
\end{eqnarray}
Since we have $p(\al) =  \frac{1}{(mn^{\al})^{1/2}}.$ Then
\begin{equation}
P[\mbox{$v_i,v_j$ are adjacent}] \leq n^{-\frac{\al}{2}}.\label{sum}
\end{equation}
Hence, for $\al \leq 2$ the above probability (\ref{sum}) is not
summable, i.e.,
\[
\sum_{n=0}^{\infty}P[\mbox{$v_i,v_j$ are adjacent}] = \infty.
\]
Since the adjacency between any pair of vertices is independent from
the adjacency between any other pair of vertices. Then by the
Borel-Cantelli Lemma, we have
\[P[\mbox{$v_i,v_j$ are adjacent} ,\qquad i.o.] = 1.\]
Therefore, for $\al \leq 2 ,$ any pair of vertices $v_i,v_j \in V,\:
i \neq j,\:i,j = 1,2,\ldots,n,$ are adjacent almost surely. Hence,
\[P[G(n,m,p(2-\ep)) \in \Theta] = 1-o(1).\]

For $\al > 2$ the above probability is (\ref{sum}) summable, i.e.,
\[
\sum_{n=0}^{\infty}P[\mbox{$v_i,v_j$ are adjacent}] < \infty.
\]
Then by the Borel-Cantelli's Lemma, we have
\[P[\mbox{$v_i,v_j$ are adjacent} ,\qquad i.o.] = o.\]
Therefore, for $\al > 2 ,$ any pair of vertices $v_i,v_j \in V$ are
adjacent only finitely many times. Hence,
\[P[G(n,m,p(2+\ep)) \in \Theta] = o(1).\]

This completes the proof. \hfill$\Box$\\

We now state strong law results for vertex degree of a typical
vertex.

Let $X = X(n,m,p)$ be the number of vertices of $V-\{v\}$ adjacent
in $G(n,m,p)$ to a vertex $v \in V,\: i,j = 1,2,\ldots,n,\:i \neq j$
i.e., $X$ be the degree of vertex $v \in V$ in $G(n,m,p).$
\begin{thm}
Let $G(n,m,p)$ be the random intersection graph with
%
$p = (mn^{\al})^{-1/2},$ where $0 < \al < 1.$ Also let $X$ be the
degree of a typical vertex $v \in V$ in $G(n,m,p).$ Then 
%
\begin{equation}
\limsup_{n \rar \infty} \frac{X}{n^{\del}} \geq a(c), \qquad a.s.,
\end{equation}
where $\del = 1-\al$ and $a(c)$ is the root in $[1,\infty)$ of
\begin{equation}
a\log\:a -a +1 = c,\label{a_c}
\end{equation}
with $a(\infty) = 1.$ \label{limsup}
\end{thm}
\textbf{Proof.} We know that the vertex degree of a vertex in graph
$G(n,m,p)$ is distributed according to Bi$(n-1,q_n).$ Hence by Lemma
\ref{lemma1}, 
we have
\begin{equation}
P[X \leq K] \leq
\exp\left((n-1)q_nH\left(\frac{(n-1)q_n}{K}\right)\right),\qquad K
\leq nq_n.
\end{equation}
From Proposition \ref{prop1}, we have for sufficiently small $p,$ we
have $q_n \sim mp^2.$ Then
%
\begin{eqnarray}
P[X \leq K] & \leq &
\exp\left((n-1)mp^2 H\left(\frac{(n-1)mp^2}{K}\right)\right)\nonumber\\
& \sim &
\exp\left(n^{1-\al} H\left(\frac{n^{1-\al}}{K}\right)\right)\nonumber\\
& = & \exp\left(-cn^{\del}\right),
\end{eqnarray}
where $H((n/mK^2)^{1/2}) = -c$ and $\del = 1-\al.$ The above
expression is summable if $K = a(c)n^{\del},$ where $a(c)$ is
increasing in $[1,\infty)$ and defined as in (\ref{a_c}). 
Then by the Borel-Cantelli Lemma we have\\
\begin{equation}
P\left[X \leq a(c)n^{\del}, \qquad i.o.\right] = 0.
\end{equation}
This implies that
\begin{equation}
\limsup_{n \rar \infty} \frac{X}{n^{\del}} \geq a(c), \qquad a.s.
\end{equation}
This completes the proof. \hfill$\Box$\\
\begin{thm}
Let $G(n,m,p)$ be the random intersection graph with
%
$p = (mn^{\al})^{-1/2},$ where $0 < \al < 1.$ Also let $X$ be the
degree of a typical vertex $v \in V$ in $G(n,m,p).$ Then 
%
\begin{equation}
\liminf_{n \rar \infty} \frac{X}{n^{\del}} \leq a(c), \qquad a.s.,
\end{equation}
where $\del = 1-\al$ and $a(c)$ is the root in $(0,1)$ of
\begin{equation}
a\log\:a -a +1 = c,\label{a_c_1}
\end{equation}
with $a(\infty) = 1.$ \label{liminf}
\end{thm}
\textbf{Proof.} We know that the vertex degree of a vertex in graph
$G(n,m,p)$ is distributed according to Bi$(n-1,q_n).$ Hence by Lemma
\ref{lemma1} 
we have
\begin{equation}
P[X \geq K] \leq
\exp\left((n-1)q_nH\left(\frac{(n-1)q_n}{K}\right)\right), \qquad
K\geq nq.
\end{equation}
From Proposition \ref{prop1}, we have for sufficiently small $p,$ we
have $q_n \sim mp^2.$ Then
%
\begin{eqnarray}
P[X \geq K] & \leq &
\exp\left((n-1)mp^2 H\left(\frac{(n-1)mp^2}{K}\right)\right)\nonumber\\
& \sim &
\exp\left(n^{1-\al} H\left(\frac{n^{1-\al}}{K}\right)\right)\nonumber\\
& = & \exp\left(-cn^{\del}\right),
\end{eqnarray}
where $H((n/mK^2)^{1/2}) = -c$ and $\del = 1-\al.$ The above
expression is summable if $K = a(c)n^{\del},$ where $a(c)$ is
decreasing in $(0,1)$ and defined as in (\ref{a_c_1}). 
Then by the Borel-Cantelli Lemma, we have\\
\begin{equation}
P\left[X \geq a(c)n^{\del}, \qquad i.o.\right] = 0.
\end{equation}
This implies that
\begin{equation}
\liminf_{n \rar \infty} \frac{X}{n^{\del}} \leq a(c), \qquad a.s.
\end{equation}
This completes the proof. \hfill$\Box$\\
%

%
\end{document}